%

\input amstex
\documentstyle{amsppt}
\magnification=\magstep1 
\pageheight{19cm}

\def\force{\Vdash}  
\def\notforce{\nVdash}
\define\k{\kappa}
\define\a{\alpha}
\redefine\l{\lambda}

\define\ra{\rightarrow}

\define\supp{\text{{\rm supp}}\,}

\define\res{\restriction}
\define\Uarr{\overrightarrow{\Cal U}}
\define\oU{o^{\Cal U}}
\define\lqot{\text{``}}
\define\rqot{\text{''}}
\define\Ult{\text{{\rm Ult}}}
\define\Gen{\text{{\rm Gen}}}

\redefine\crit{\text{{\rm crit}}\,}
\define\Reg{\text{{\rm Reg}}}
\define\Sing{\text{{\rm Sing}}}
\define\Club{\text{{\rm Club}}}

\define\Tr{\text{{\rm Tr}}\,}
\define\CU{\text{{\rm CU}}}

\topmatter
\title
Consistency Strength of the Axiom of Full Reflection at Large Cardinals
\endtitle

\author Moti Gitik and Jiri Witzany \endauthor


\affil  Tel Aviv University \\
 The Pennsylvania State University and Charles University (Prague)
             \endaffil

\address    School of Mathematical Sciences, Sackler Faculty 
of Exact Sciences, Tel Aviv University, Ramat-Aviv, 69978 Israel
\endaddress
\email  gitik\@math.tau.ac.il \endemail

\address
Department of Mathematics, The Pennsylvania State University,
University Park, PA 16802
\endaddress
\email  witzany\@math.psu.edu \endemail



\abstract We prove that the Axiom of Full Reflection at 
a measurable cardinal is equiconsistent with the existence
of a measurable cardinal. We generalize the result also
to larger cardinals as strong or supercompact.
    \endabstract 

\date April 20, 1994 \enddate

\endtopmatter

\def\lheadline{\folio\hfil {\eightpoint MOTI GITIK AND JIRI WITZANY}\hfil}
\def\rheadline{\hfil{\eightpoint ON REFLECTION OF STATIONARY SETS} \hfil\folio}
\headline={\ifodd\pageno\rheadline \else\lheadline\fi}

\document
\baselineskip=15pt

\head {\bf 1. Introduction} \endhead

It has been proved in [JS90] that the Axiom of Full Reflection at an 
$n$-Mahlo cardinal is equiconsistent with a $\Pi^1_n$-indescribable
cardinal and in [JW94] that  consistency of the Axiom of Full
Reflection at a measurable cardinal follows from consistency
of a coherent sequence of measures with a repeat point.
It has been conjectured in [JW94] that the two principles are
actually equiconsistent. However we prove that Full Reflection at 
a measurable cardinal can be obtained surprisingly from
only one measure. Furthermore the method also
generalizes to larger cardinals as strong or supercompact.
Hence we can conclude that the Axiom of Full Reflection at large cardinals
weaker than measurable, e.g. as $n$-Mahlo, does push the
 consistency strength up,
but does not push  the consistency strength up at a measurable
or larger cardinals.

To state the main theorem let us review the basic definitions and facts.
If $S$ is a stationary subset of a regular uncountable cardinal $\k$ then
{\it the trace of} $S$ is the set

$$ \Tr (S) = \{ \a < \k ; \; \; S \cap \a \; \text{is stationary
in}  \; \a \} $$
and we say that $S$ {\it reflects at} $\a\in \Tr(S)$.  If $S$ and $T$ are
both stationary,
we define
$$ S < T \;\; \text{if for almost all} \; \a \in T, \; \; \a \in
\Tr(S) $$
and say that $S$ {\it reflects fully} in $T$.
(Throughout the paper, ``for almost all" means ``except for a
nonstationary set of
points"). It can be proved that this relation is a well-founded
partial ordering (see [JW94] or [J84]).
 {\it The order} $o(S)$ of a stationary set of regular cardinals is defined
as the rank of S in the relation $<$:
$$ o(S)=\sup \{ o(T)+1 ; \; T \subseteq \Reg(\k ) \text{ is stationary and }  
T<S \} .$$
For a stationary set $T$ such that $T\cap \Sing(\k )$ is stationary 
define $o(T)=-1$. {\it The order of } $\k$ is then defined as
$$ o(\k) = \sup\{ o(S)+1 \; ;  S \subseteq \k \; \text{is stationary} \}.$$
Note that the order
$o(\k)$ provides a natural generalization of the Mahlo hierarchy: $\k$ is
exactly $o(\k)$-Mahlo if $o(\k) < \k^{+}$ and  greatly Mahlo
if $o(\k) \geq \k^{+}$.

We say that a stationary set $S$ {\it reflects fully} at regular cardinals 
if for
any stationary set $T$ of regular cardinals $o(S)<o(T)$ implies $S<T$.

\proclaim{Axiom of Full Reflection at $\k$}  Every stationary subset of $\k$ 
reflects
fully at regular cardinals. \endproclaim

Notice that the axiom presents in a sense the maximal possible amount
of reflection of stationary subsets of $\k$ at regular cardinals.

Now we are able to state the main theorem:

\proclaim{Theorem} Let $\phi(\k)$ be one of the following principles:

(i) $\k$ is measurable,

(ii) the Mitchell order of $\k$ is $\k^{++},$

(iii) $\k$ is $n$-strong,

(iv) $\k$ is strong,

(v) $\k$ is $\k^{+n}$-supercompact,

(vi) $\k$ is supercompact.

Assume that $V$ satisfies GCH and $\phi(\k),$ then there is a model
where GCH, the Axiom of Full Reflection at $\k,$ and $\phi(\k)$ hold.
\endproclaim

The case (ii) has been actually proved in [JW94]: it has been proved
in the paper that if $\Uarr$ is a coherent sequence of measures
then there is a forcing notion $P_{\k+1}$ that preserves any 
repeat point of $\Uarr$ on $\k.$ If $\oU(\k)=\k^{++}$ then there are
$\k^{++}$ repeat points on $\k$ and it is not difficult to see that
the Mitchell order of $\k$ is $\k^{++}$ in the generic
extension by $P_{\k+1}.$ Thus we will work only on cases (i) and (iii)-(vi).

\head {\bf 2. Proof of the theorem} \endhead

The proof should be self-contained, however a knowledge 
of [JW94] is helpful.

Assume that $V$ satisfies GCH and $j: V \ra M$
is an elementary embedding such that $\crit(j)=\k$ and
$V \cap \;^\k M \subseteq M .$ We will define a forcing
$P_{\k+1}$ that will work in all cases (i),(iii)-(vi).
$P_{\k+1}$ will be an Easton support iteration of 
$\langle Q_\l; \l\leq \k \rangle ;$
$Q_\l$ will be nontrivial only for $\l$ Mahlo, and in that
case it will be an iteration (defined in $V(P_\l)$)
of length $\l^+$ with $<\l$ support of forcing notions
shooting clubs through certain sets $X\subseteq \l$
always with the property that $X\supseteq \Sing(\l) .$
This will be guarantee $Q_\l$ to be essentially 
$<\l$-closed (i.e. it will have a $<\l$-closed dense subset). Consequently $Q_\l$ will be $\l^+$-c.c.,
$P_\l$ will be $\l$-c.c., and the factor iteration 
$P_{\l+1,\k+1}$ above $\l$ will be essentially $\l$-closed.
Therefore $P_{\k+1}$ will preserve cardinals, cofinalities, and GCH.

Consider an iteration $Q$ of 
$\langle \CU(\dot X_\a); \a<l(Q)\rangle$ with $<\l$ support,
where $\CU(\dot X_\a)$ denotes the forcing shooting a club in 
$V(P_\l\ast Q\res\a)$ through a subset $\dot X_\a$ of $\l$
containing $\Sing(\l).$ In that case we say that $Q$ is an iteration of order
0. Since $Q\res\a$ is essentially $<\l$-closed, conditions in
$\CU(\dot X_\a)$ can be taken in $V(P_\l)$ rather than in 
$V(P_\l\ast Q\res\a).$ So $Q$ can be considered to be a set of
sequences of closed bounded subsets of $\l$ in $V(P_\l).$
Since $P_\l$ is $\l$-c.c. there is an appropriate $P_\l$-name for $Q$
of cardinality
$\l$ if $l(Q)<\l^+,$ and of cardinality $\l^+$
if $l(Q)=\l^+ .$ Let $\tilde Q$ be another iteration of
$\langle \CU(\dot Y_\gamma); \gamma<l(\tilde Q) \rangle$
of order 0. We say that $Q$ is an subiteration of $\tilde Q$
if there is a 1-1 function $\pi:l(Q)\ra l(\tilde Q):\a\mapsto \gamma_\a$
inducing an embedding of $Q$ into $\tilde Q$
such that
$\dot X_\a$ is an equivalent name to $\dot Y_{\gamma_\a}$ with
respect to the induced embedding of $Q\res\a$ into $\tilde Q.$
Notice that the sequence $\langle \gamma_\a; \a<l(Q) \rangle$
 does not have to be increasing.
Any $Q$-name can be considered to be  a $\tilde Q$-name via
the induced embedding; $\tilde Q$ is actually isomorphic
to an iteration of order 0 in the form $Q\ast R.$

We will need to estimate (in $V(P_\l)$) the number of iterations
of order 0 and length $<\l^+.$ Each such iteration is a set
of sequences with $<\l$ support of bounded subsets of $\l.$
Therefore it is easy to see that the number is $\leq 2^\l=\l^+.$

For any iteration $Q$ of order $\delta +1$ we will define certain filters
$F^Q_{\l,\delta}$ on $\l$ in $V(P_\l\ast Q).$
Simultaneously by induction on $\beta$ and $l(Q)$ we define $Q$ to
be an iteration of order $\beta$ if it is an iteration of
$\langle CU(\dot X_\a); \a<l(Q) \rangle$ 
with $<\l$-support such that $l(Q)<\l^+$ and for all $\a<l(Q)$:
$$P_\l\ast Q\res\a \force \lqot \Sing(\l)\subseteq\dot X_\a
\text{ and } \dot X_\a \in F^{Q\res\a}_{\l,\delta} 
\text{ for all } \delta<\beta .\rqot$$
Let us call such an assignment $Q\mapsto F^Q_{\l,\delta}$
a filter system $F_{\l,\delta} .$
$F_{\l,\delta}$ will be defined for all $\delta<\Theta(\l)$
where $\Theta(\l)$ will be specified later.
The filter systems will have among others the property that
$F^{Q\res\a}_{\l,\delta}\subseteq F^Q_{\lambda,\delta}.$

$Q_\l$ is then defined in $V(P_\l)$ to be an iteration
of length $\l^+$ such that for all $\a<\l^+$
$Q_\l\res\a$ is an iteration of order $\Theta(\l),$ and all
potential names for subsets of $\l$ are used cofinally many times.

It remains to find the filter systems $F_{\l,\delta}$
(working in $V(P_\l)$). We require that for any iteration $Q$
of order $\delta +1$ the following is satisfied:
\vskip5pt
\noindent (i) If $Q'$ is an subiteration of $Q$ then
$$F^{Q'}_{\l,\delta} = F^Q_{\l,\delta}\cap V(P_\l\ast Q'),$$
(ii) $P_\l\ast Q \force \lqot F^Q_{\l,\delta}\supseteq \Club(\l)
\text{ is a proper filter,}$\newline
$$\forall S\subseteq\Sing(\l)\text{ stationary: }\Tr(S)\in F^Q_{\l,\delta},$$
$$\forall S\subseteq\l:(\exists\gamma<\delta:S\text{ is }
F^Q_{\l,\gamma}\text{-positive}) \Rightarrow \Tr(S)\in F^Q_{\l,\delta},\rqot$$
(iii) $P_\l\ast Q\force \lqot \forall S\subseteq \Reg(\l):
(\forall\gamma<\delta: S\text{ is }F^Q_{\l,\gamma}\text{-thin})
\Rightarrow \k\setminus\Tr(S) \in F^Q_{\l,\delta}.\rqot$
\vskip5pt
\noindent Moreover we require that
\vskip5pt
\noindent (iv) there is an iteration $Q$ of order $\delta+1,$ a 
$P_\l\ast Q$-name $\dot X$ for a subset of $\l$ and $p\ast q\in P_\l\ast Q$
so that
$$p\ast q\force_{P_\l\ast Q} \lqot \dot X \text{ is $F^Q_{\l,\gamma}$-thin
for all $\gamma<\delta$,'' but}$$
$$p\ast q\force_{P_\l\ast Q} \lqot \dot X \text{ is 
$F^Q_{\l,\delta}$-positive.''}$$
\vskip5pt
By in induction on $\delta$ choose a filter system $F_{\l,\delta}$
as long as there is such a filter system  with properties (i)--(iv).
Since  the number of iterations $Q$ of length
$<\l^+$ with $<\l$-support shooting closed unbounded subsets of $\l$
is $\leq \l^+$
and since $F_{\l,\delta}$ is by (iv) different from all
$F_{\l,\gamma}$ ($\gamma<\delta$), this process must eventually stop
after a number of steps $\Theta(\l)<\l^{++}.$

Apply this process by induction on all $\l<\k$ defining
an iteration $P_\k$ below $\k.$ Put $P_{\k+1}=(jP_\k)\res(\k+1).$
Note that $P_{\k+1}=P_\k\ast Q_\k$ where $Q_\k$ is an iteration
of length $\k^+$ with $<\k$-support, given by certain filter systems 
 $F_{\k,\delta}$ ($\delta<\Theta=\Theta(\k)$). 

We claim that
$$V(P_{\k+1})\models\text{``Full Reflection at $\k$''}$$
and the embedding $j$ can be in many cases lifted onto $V(P_{\k+1}).$

Let us define $F_j$ in $V$ similarly as in [JW94] to be a $\Theta$-th
filter system on $\k$:

By induction on $l(Q)$ say that $Q,$ an iteration
of $\langle CU(\dot X_\a); \a<l(Q) \rangle,$ is an iteration of order
$\Theta +1$ w.r.t. $F_j$ if it is an iteration of order $\Theta$ and for 
all $\a<l(Q)$
$$P_\k \ast Q\res\a \force \lqot \dot X_\a \in F^{Q\res\a}_j .\rqot$$
If $Q$ is an iteration of order $\Theta +1$ w.r.t. $F_j,$
 $\dot X$ a $P_\k\ast Q$-name,
$p\ast q\in P_\k\ast Q,$ we define 
$p\ast q\force\lqot \dot X \in F^Q_j \rqot$ if
$$p\force_{jP_\k}\lqot \forall H\in \Gen_j(Q,G^*): q \in H \Rightarrow
[H]^j\force_{jQ} \k \in j\dot X .\rqot \tag1$$
Here $\Gen_j(Q,G^*)$ is defined as follows:
let $G^*$ be a $jP_\k$-generic filter over $V,$ $G=G^*\res P_\k.$
Then $Q$ is obviously an subiteration of $Q_\k$ which gives
a  filter $H$ from $G^*\res Q_\k$ that is $Q$-generic over $V[G].$
$\Gen_j(Q,G^*)$ denotes the set of all filters $H$ obtained in this way.
We can easily find many $H\in \Gen_j(Q,G^*)$ such that $q\in H$:
since $Q_\k$ is an iteration of order $\Theta$ such that all potential
names are used cofinally many times we can find a sequence
of ordinals $\langle \gamma_\a; \a<l(Q) \rangle$
inducing a subiteration embedding of $Q$ into $Q_\k$ such that
all $\gamma_\a$'s are above any given $\beta<\k^+$;
hence by a density argument there is $r\in G^*$ and such a sequence
$\langle \gamma_\a; \a<l(Q) \rangle$ with the property that
$r\res \langle \gamma_\a; \a<l(Q) \rangle = q.$

Represent an $H\in \Gen_j(Q,G^*)$ as 
$\langle C_\beta;\beta<l(Q) \rangle$ where $C_\beta$'s
are the generic closed unbounded subsets of $\k.$
$[H]^j$ is a sequence of length $j(l(Q))$ defined as follows
$$[H]^j(\gamma)=\cases C_\beta\cup\{\k\}\text{ if }j(\beta)=\gamma, \\
                       \emptyset \text{ otherwise.}
                \endcases $$
To prove that $[H]^j \in jQ$ all we need is to check inductively
that $[H]^j\res j(\beta)\force_{j(Q\res\beta)} \lqot \k\in j\dot X_\beta.\rqot$
But this immediately follows from the assumption
$P_\k\ast Q\res\beta \force \lqot \dot X_\beta \in F^{Q\res\beta}_j .\rqot$

\proclaim{Lemma 1} The filter system $F_{\k,\Theta}=F_j$ satisfies (i)--(iii)
with
$\delta=\Theta .$
\endproclaim

\demo{Proof} (i) Let $Q,Q'$ be two iterations of order $\Theta+1;$ 
 assume $\pi$ embeds $Q$ into $Q'$ via 
$\langle \a_\delta;\delta<l(Q) \rangle$ as an subiteration.
Let $\dot X$ be a $P_\l\ast Q$-name for a subset of $\l.$

Suppose $ p\ast q\in P_\l\ast Q,$
$ p\ast q\force_{P_\l\ast Q}\lqot \dot X\in F_j^Q .\rqot $
We want to prove that 
$$p\ast \pi(q)\force _{P_\l\ast Q'}\lqot \pi(\dot X)\in F_j^{Q'}  .\rqot $$
Let $G^*\ni p$ be $jP_\l$-generic over $V,$ 
$H'\in \Gen_j(Q',G^*),$ $H'\ni \pi(q) .$
Then the embedding of $Q'$ into $(jP_\l)^\l$ induces via $\pi$ 
an embedding of $Q$ into $(jP_\l)^\l$ giving $H\in \Gen_j(Q,G^*)$
such that $q\in H .$
Moreover $j\pi$ embeds $jQ$ into $jQ'$ by elementarity, and
$(j\pi)([H]^j)\geq[H']^j .$
Since $[H]^j\force_{jQ}\lqot\check \l\in j \dot X \rqot $
it follows that 
$[H']^j\force_{jQ'}\lqot\check \l\in j(\pi \dot X ) .\rqot $

Now suppose $p\ast q'\in P_\l\ast Q',$ 
$p\ast q'\force _{P_\l\ast Q'}\lqot \pi(\dot X)\in F_j^{Q'} .\rqot $
Let $q\in Q$ be such that $\pi(q)$ agrees with $q'$ on the set
$\{\a_\delta;\delta<l(Q) \} .$
We claim that $ p\ast q\force_{P_\l\ast Q}\lqot \dot X\in F_j^Q .\rqot $
Let $G^*\ni p$ be $jP_\lambda$-generic over $V,$
$H\in \Gen_j(Q,G^*),$ and $q\in H.$
We need to prove $[H]^j\force_{jQ}\lqot\check \l\in j \dot X .\rqot $
Suppose it is not true, then there is $\tilde q \leq[H]^j$
such that $\tilde q\force_{jQ}\lqot\check \l\notin j \dot X .\rqot $
Express $Q'$ as $Q\ast R,$ and as above  find
a subiteration embedding of $Q'$ into $(jP_\l)^\l$ that extends the embedding
of $Q$, giving $H'\in \Gen_j(Q',G^*)$ such that 
$H'\res Q = H,$ and $q'\in H'.$
In other words if $\pi_1 :l(Q)\ra \l^+$ embeds $Q$ into $(jP_\l)^\l$
then we obtain $\pi_2:l(Q')\ra \l^+$ embedding $Q'$ into $(jP_\l)^\l$
such that $\pi_1(\delta)=\pi_2(\a_\delta)$ for $\delta < l(Q) .$ 
Now $j\pi$ embeds $jQ$ into $jQ'$ via 
$j\langle \a_\delta;\delta<l(Q) \rangle ,$ thus
$(j\pi)(\tilde q)\in jQ'$ and 
$$\supp((j\pi)(\tilde q))\subseteq j(\{\a_\delta;\delta<l(Q)\}).$$
Moreover $\supp([H']^j)=j''l(Q'),$ if $\a<l(Q')$
then either $\a\in \{\a_\delta;\delta<l(Q)\},$ and then
$(j\pi)(\tilde q)(j\a)$ extends $[H']^j(j\a),$
or $\a\notin \{\a_\delta;\delta<l(Q)\},$ then
$j(\a)\notin \supp((j\pi)(\tilde q)).$
Consequently $(j\pi)(\tilde q)$ and $[H']^j$ are compatible.
But $[H']^j\force_{jQ'}\lqot \check \l\in j(\pi \dot X),\rqot $
while $(j\pi)(\tilde q)\force_{jQ'}\lqot \check \l\notin j(\pi \dot X)\rqot $
 - a contradiction.

(ii) Each $F^Q_j$ is obviously proper and contains $\Club(\k).$
Let $P_\k\ast Q \force \lqot \dot S \subseteq \k
\text{ is $F^Q_{\k,\gamma}$-positive}\rqot$ for some $\gamma<\Theta$
(or $\dot S\subseteq\Sing(\k)$ is stationary).
We wish to prove that $P_\k\ast Q \force \lqot \Tr(\dot S)\in F^Q_j .\rqot$
Assume towards a contradiction that $G^*$ is $jP_\k$-generic over $V,$
$H\in \Gen_j(Q,G^*),$ and 
$[H]^j \notforce _{jQ} \lqot \k \in \Tr(\dot S) .\rqot$
So there is $H^* \ni [H]^j$ $jQ$-generic over $V[G^*]$ so that
$$V[G^*\ast H^*]\models \lqot S\text{ is nonstationary.}\rqot$$
Since $(jP_\k)_{\k+1,j\k} \ast jQ$ is essentially $\k$-closed
and $Q_\k$ is $\k^+$-c.c. there is a
sufficiently large $\a<\k^+,$ such that if $G=G^*\res P_\k,$
$\tilde H =G^*\res (Q_\k\res \a)$ then
$$V[G\ast \tilde H]\models \lqot S \text{ is nonstationary,}\rqot$$
which is a contradiction with (i) as $V[G\ast H]\models$``$S$ is 
$F^H_{\k,\gamma}$-positive'' and $Q$ is an subiteration of $Q_\k\res\a $
giving $H$ from $\tilde H$ (provided $\a$ is large enough).

(iii) Assume that
$$P_\k\ast Q\force \lqot \dot S \subseteq\Reg(\k)\text{ and }
\forall\gamma<\Theta : \dot S \text{ is }F^Q_{\k,\gamma}
\text{-thin}.\rqot$$
We want to prove that 
$P_\k\ast Q\force \lqot \k\setminus\Tr(\dot S)\in F^Q_j .\rqot$
Assume $G^*$ is $jP_\k$-generic, $H\in \Gen_j(Q,G^*),$
$H^*\ni [H]^j$ $jQ$-generic over $V[G^*]$ and
$V[G^*\ast H^*]\models\lqot \k \notin j(\k\setminus\Tr(S)),\rqot$
i.e. $V[G\ast \tilde H]\models$``$S$ is stationary''
where $\tilde H = (G^*)\res Q_\k.$
But a club have been shot through $\k\setminus S$ in the iteration 
$Q_\k$ - a contradiction.
\qed
\enddemo

\proclaim{Lemma 2} Let $Q_o$ be an iteration of 
$\langle CU(\dot X_\a); \a<l(Q) \rangle $ of
order $\Theta,$
$\dot X$ a $P_\k\ast Q_o$-name for a subset of $\k,$ 
$p\ast q \in P_\k \ast Q_o.$
Then $Q_o$ is an iteration of order $\Theta+1$ w.r.t. $F_j,$
and moreover
if $p\ast q\force$``$\dot X$ is $F^{Q_o}_{\k,\gamma}$-thin for 
all $\gamma<\Theta$''
then $p\ast q \force \lqot \dot X$ is $F^{Q_o}_j$-thin.''
\endproclaim

\demo{Proof} Assume towards a contradiction that
$p\ast q\force \lqot\dot X$ is $F^{Q_o}_j$-positive.''
Then we claim that the construction of filter systems $F_{\k,\gamma}$
in $M=\Ult(V,U)$ could not stop at $\Theta .$
$F_j$ cannot be constructed in $M,$
however we can construct its approximation.

Firstly define $\tilde F_{\k , \Theta}$ as follows:

Let $\tilde F^\emptyset_{\k , \Theta}$ ($Q=\emptyset$) be generated in $V(P_\k)$
by all sets that should be there by (ii) and (iii), and
by $\dot X_0.$ Note that $\dot X_\a$ is forced to be in $F^{Q_o\res \a}_j$
for all $\a<l(Q)$ by the induction hypothesis. Hence 
$\tilde F^\emptyset_{\k , \Theta}
\subseteq F^\emptyset _j$ verifying that 
$\tilde F^\emptyset_{\k , \Theta}$ is a 
proper filter.
Similarly define $\tilde F_{\k , \Theta}^Q$ for iterations $Q$
of order $\Theta+1$ w.r.t. previously defined 
$\tilde F^{Q\res\a}_{\k , \Theta}.$
We also have to make sure that 
$\dot X_\a \in  \tilde F^{Q_o\res\a}_{\k , \Theta}$
for all $\a<l(Q_o).$
This filter system satisfies (ii) and (iii), clearly
$\tilde F^{Q'}_{\k , \Theta}\subseteq \tilde F^{Q}_{\k , \Theta}$
if  $Q'$ is an subiteration of $Q,$ however (i) does not have to hold.
To achieve that define
$$F^\emptyset_{\k,\Theta}=\bigcup \{ \tilde F^{Q}_{\k , \Theta}\cap V(P_\k);
  \text{ $Q$ is an iteration of order 
        $\Theta+1$ w.r.t. $\tilde F_\Theta$ } \}.$$
Then for $Q$ an iteration of order $\Theta+1$ w.r.t. previously defined
$F^{Q\res\a}_{\k,\Theta}$'s by induction on $l(Q)$ define
$$F^Q_{\k,\Theta}=\bigcup \{ \tilde F^{Q'}_{\k , \Theta}\cap V(P_\k\ast Q);
  \text{ $Q'$ is an iteration of order $\Theta+1$ w.r.t. 
$\tilde F_{{\k , \Theta}}$ }$$
$$\text{such that $Q$ is an subiteration of $Q'$}\}.$$
It is not difficult to see that such $Q'$ exists.
We have constructed a filter system $F_{\k,\Theta}$ in $M$ that
satisfies (i)--(iii).
Moreover $Q_o$ is an iteration of order $\Theta+1$ w.r.t. $F_{\k,\Theta},$
and so (iv) holds for the $\dot X,$ $p\ast q$ from the assumption
of the lemma - a contradiction.
\qed
\enddemo

Let $G\ast H$ be $P_\k\ast Q_\k$-generic over $V.$

\proclaim{Lemma 3} $V[G\ast H]\models$``Full Reflection holds up to $\k.$''
\endproclaim

\demo{Proof} For $\gamma<\Theta$ define 
$F_{\k,\gamma}^H=\bigcup_{\a<\k^+} F^{H\res\a}_{\k,\gamma}.$
We know that $F_{\k,\gamma}^H\supseteq \Club(\k)$ is proper. By (i)
if $S\in V[G\ast H\res\a]$ is $F_{\k,\gamma}^{H\res\a}$-positive
then it is $F_{\k,\gamma}^H$-positive. Moreover by the construction
$S\subseteq \Reg(\k)$ is stationary iff $S$ is $F_{\k,\gamma}^H$-positive
for some $\gamma<\Theta$ iff $S$ is $F_{\k,\gamma}^{H\res\a}$-positive
whenever $S\in V[G\ast H\res\a].$ Let us firstly prove 
that $V[G\ast H]\models \lqot S<\Reg(\k) \rqot$
for $S\subseteq\Sing(\k)$ stationary in $V[G\ast H].$ Let 
$S\in V[G\ast H\res\a]$ then $S$ is also stationary in this model,
and so by (ii) $\Tr(S)\in F^{H\res\a}_{\k,\gamma}$ for all $\gamma<\Theta,$
consequently a club has been shot through 
$\Sing(\k)\cup \Tr(S).$

Now let $S\subseteq\Reg(\k)$ be stationary, denote $\gamma_S$
to be the least $\gamma$ such that $S$ is $F_{\k,\gamma}^H$-positive.
The following claim completes the proof of Full Reflection at $\k$
in $V[G\ast H]$ (the proof for $\l<\k$ is identical).

\proclaim{Claim} Let $S,T \subseteq \Reg(\k)$ be two stationary
sets. Then $\gamma_S < \gamma_T$ iff $S<T.$ Consequently $\gamma_S=\gamma_T$
iff $o(S)=o(T).$
\endproclaim

\demo{Proof} Let $S,T \in V[G\ast H\res\a],$ $\gamma_S<\gamma_T.$
Then $S$ is $F_{\k,\gamma_S}^{H\res\a}$-positive, and so by (ii)
$\Tr(S)\in F_{\k,\delta}^{H\res\a}$ for all $\delta>\gamma_S.$
Thus $T\setminus\Tr(S)$ is $F_{\k,\delta}^{H\res\a}$-thin
for all $\delta < \Theta,$ so a club has been shot through
$\k\setminus(T\setminus\Tr(S)),$
which means that $T\setminus\Tr(S)$ is nonstationary in $V[G\ast H],$
i.e. $S<T.$

On the other hand assume that $S<T,$ then necessarily 
$\gamma_S \leq \gamma_T .$ By the definition of $\gamma_S$
the set $S$ is $F_{\k,\delta}^{H\res\a}$-thin
for all $\delta<\gamma_S,$ and so by (iii)
$\Tr(S)$ is $F^{H\res\a}_{\k,\gamma_S}$-thin.
Since
$T\setminus\Tr(S)$ is nonstationary in 
$V[G\ast H],$ it must be $F_{\k,\gamma_S}^{H\res\a}$-thin.
Thus $T=(T\setminus\Tr(S))\cup\Tr(S)$ is $F_{\k,\gamma_S}^{H\res\a}$-thin
proving $\gamma_S <\gamma_T .$

Finally if $\gamma_S=\gamma_T$ and say $o(S)<o(T)$ then
there must be $S'<T$ such that $o(S)=o(S').$
By the fact proven above $\gamma_{S'} < \gamma_T =\gamma_S,$
and so $S'<S$ - a contradiction.\newline
\qed$\;$Claim, Lemma 3
\enddemo
\enddemo

Finally we need to prove that $P_{\k+1}$ preserves large cardinal
properties of $\k.$ Let us firstly consider measurability and
supercompactness of $\k.$

\proclaim{Lemma 4} Let $\l\geq\k$ be a cardinal such that

(i) $V\cap\;^\l M \subseteq M ,$

(ii) $\l^+ < j(\k) < j(\k^+) < \l^{++} ,$

(iii) there is no Mahlo cardinal between $\k$ and $\l+1.$ 

\noindent Then the embedding $j:V\ra M$ can be extended to
$j^{**}:V[G\ast H] \ra M[G^* \ast H^*]$ in $V[G\ast H]$ so that
$V[G\ast H] \cap\;^\l M[G^*\ast H^*]\subseteq M[G^*\ast H^*].$
\endproclaim

\demo{Proof} By the definition of $P_{\k+1}$ the forcing
$jP_{\k+1}$ factors as $P_{\k+1}\ast R_o\ast j(Q_\k) .$
So all we need is to find an $R_o\ast j(Q_\k)$-generic
filter $H_o\ast H^*$ over $M[G\ast H]$ so that $p\ast q\in G\ast H$
implies $j(p\ast q )\in G\ast H\ast H_o \ast H^* .$
The factor iteration $R_o=(jP_{\k+1})_{\k+1, j\k}$ starts with
a nontrivial forcing at the first Mahlo cardinal in $M$ above
$\k$ which must be above $\l.$ Consequently $R_o$ is essentially
$\l$-closed in $M[G\ast H]$ as well as in $V[G\ast H].$
Let $D$ be a $\l$-closed dense subset of $R_o .$
The number of dense subsets of $D$ in $M[G\ast H]$ is $j(\k^+)$
and the cardinality of $j(\k^+)$ in $V$ is just $\l^+ .$
Thus we have only $\l^+$ dense subsets of a forcing that is $\l$-closed
in $V[G\ast H],$ and so it is easy to construct
$H_o\in V[G\ast H]$ that is $R_o$-generic over $M[G\ast H].$
Obviously $p\in G$ implies $j(p)\in G^*=G\ast H \ast H_o ,$
thus $j$ extends to $j^* : V[G] \ra M[G^*]$ in $V[G\ast H].$
It immediately follows from the $\k$-c.c. of $P_\k$ that
$V[G] \cap \; ^\l M[G^*]\subseteq M[G^*] .$
Next we need to find a filter $H^* \in V[G\ast H]$ that is  
$j^*(Q_\k)$-generic over $M[G^*],$ and such that
$[H\res\a]^j \in H^*$ for all $\a<\k^+ .$

It is easy to see that the number of antichains of $Q_\k$
(in $V[G]$) is only $\k^+$: if $A\subseteq Q_\k$ is an antichain,
then $|A|\leq\k,$ which implies that there is an $\a<\k^+$ such that
$A\subseteq Q_\k\res\a,$ the number of subsets of $Q_\k\res\a$ is only
$\k^+ .$
By elementarity $M[G^*]\models$``the number of antichains in $j^*(Q_\k)$
is $j(\k^+)$''. Moreover $M[G^*]\models$``$j^*(Q_\k)$ is essentially
$\l$-closed''. Let $D$ be a $\l$-closed dense open subset
of $j^*(Q_\k),$ put
$$\Cal D =\{A\in M[G^*]; A\subseteq D \text{ is an antichain} \} .$$
Then $V[G\ast H]\models \lqot D$ is $\l$-closed, 
$|\Cal D|=|j(\k^+)|=\l^+ .\rqot $
Now we have to distinguish two cases: if $\l\geq\k^+$ then 
$[H]^j=\cup_{\a<\k^+}[H\res\a]^j$ is a good master condition
in $j^*(Q_\k),$ and we can easily build up $H^*\in V[G\ast H]$
$j^*(Q_\k)$-generic over $M[G^*]$ such that $[H]^j\in H^* .$
If $\l=\k$ then we have to be more careful. Let
$\langle A_\a;\a<\k^+\rangle$ be an enumeration of $\Cal D$ in which
each element of $\Cal D$ occurs cofinally many times. Construct a
descending sequence of conditions 
$\langle q_\a; \a<\k^+ \rangle \subseteq D$ with the following
properties

(i) $q_\a \in j^*(Q_\k\res\a) ,$

(ii) $q_\a \leq [H\res\a]^j,$

(iii) if $A_\a\subseteq j^*(Q_\k\res\a)$ then $q_\a$ strengthens a 
condition in $A_\a .$

\noindent The sequence $\langle q_\a;\a<\k^+\rangle$ generates a
$j^* Q_\k$-generic filter $H^*\in V[G\ast H]$ over $M[G^*]$
such that each  $[H\res\a]^j$ is in $H^* .$

Since $P_{\k+1}$ is $\k^+$-c.c. each $P_{\k+1}$-name for
a $\l$-sequence of ordinals in $V$ is already in $M.$
Hence $V[G\ast H]\cap\; ^\l M[G^* \ast H^*] \subseteq M[G^* \ast H^*] .$
\qed
\enddemo

By the lemma if $\k$ is measurable, or $\l$-supercompact
with no Mahlo cardinal between $\k$ and $\l+1,$ and if 
$P_{\k+1}$ is constructed using a corresponding elementary embedding $j,$
then the forcing preserves measurability, or $\l$-supercompactness of $\k.$

Now suppose $\k$ is supercompact. We can assume without loss
of generality that there is no inaccessible cardinal above $\k,$
cutting off the universe if there is any.
For each $\l> \k$ there is a $\l$-supercompact embedding $j$
given by an ultrafilter on $\Cal P_\k(\l).$
Assign to $\l$ a forcing $P^\l_{\k+1}$ constructed from $j$
as above. It is easy to estimate the number of possible forcings $P_{\k+1}$
to be $\leq \k^{++}.$ Consequently there is a proper class of
$\l$'s with the same  $P_{\k+1}=P^\l_{\k+1} .$ This $P_{\k+1}$
preserves $\l$-supercompactness of $\k$ for any of those $\l$'s, and so
supercompactness of $\k .$

Let us turn our attention to strong cardinals. The following
is essentially the idea how to modify the construction above.

\proclaim {Lemma 5} Let $j:V\ra M$ be given by a $(\k,\l)$-extender:
$\crit(j)=\k,$ $V\cap\; ^\k M\subseteq M,$
$M=\{ (jf)(a);  a\in [\l]^{<\omega} , f\in \;^{[\k]^{|a|}} V \}.$
Moreover assume that $P$ is a notion of forcing such that
$M\models$``$|P|\leq j(\k^+),$ $P$ has $j(\k^+)$-c.c., and
$P$  is $\l$-closed.''
Then there is $G\in V$ $P$-generic over $M.$
\endproclaim

\demo{Proof} (J. Zapletal) 
We can assume that $P\subseteq j(\k^+) .$
Let $\langle f_\a; \a<\k^+ \rangle$ be an enumeration of all
functions $\k\ra [\k^+]^\k .$
Construct a sequence $\langle p_\a; \a< \k^+ \rangle$ of conditions in $P$
as follows: Put $p_0=1.$ For limit $\a$ get a lower bound of
$\langle p_\delta; \delta<\a \rangle$ using closedness of $M$
and $P.$ For $\a=\beta+1$ put
$X = \{ (jf_\beta)(a); a\in [\l]^{<\omega},$ $(jf_\beta)(a)\subseteq P$
is a maximal antichain$\}.$
$X$ is a set in $M$ of cardinality $\leq \l,$ hence we can find
$p_{\beta+1}<p_\beta$ that meets all of those maximal antichains
using closedness of $P$ in $M.$

By the chain condition the filter $G$ generated by
$\langle p_\a;\a<\k^+ \rangle$ is $P$-generic over $M.$
\qed
\enddemo

Let $j: V \ra M$ be $\gamma$-strong, i.e. $\crit(j)=\k,$
$V_{\k+\gamma}\subseteq M,$ $\gamma < j(\k) .$
It is a standard fact on extenders (see [Ka93])
that we can assume 
$$M=\{(jf)(a); a\in [\l]^{<\omega}, f\in \;^{[\k]^{|a|}}V \} ,$$
where $\l=|V_{\k+\gamma}|^{+M} < j(\k) .$

Assume there is no Mahlo cardinal between $\k$ and $\l+1 .$
Let $P_{\k+1}$ be constructed from $j,$
$jP_{\k+1} = P_{\k+1} \ast R_o \ast (jQ_\k) ,$
$G\ast H$ $P_{\k+1}$-generic over $V.$
To construct $H_o \in V[G\ast H]$ $R_o$-generic over $M[G\ast H]$
consider an enumeration 
$\langle f_\a; \a<\k^+ \rangle$ of all functions in $V$ from
$\k$ to $[P_\k]^\k .$
Construct a descending chain 
$\langle p_\a;\a<\k^+ \rangle \subseteq R_o$ similarly as in the
proof of lemma 5 so that $p_\a$ meets any maximal
antichain $\subseteq R_o$ of the form
$(jf_\a)(a)/G\ast H$ ($a\in [\l]^{<\omega}$).
We only have to observe that $R_o$ is $\k$-closed in $V[G\ast H]$
and $\l$-closed in $M[G\ast H].$
The sequence $\langle p_\a; \a<\k^+ \rangle$ generates
a filter $H_o\subseteq R$ generic over $M[G\ast H].$
Now $j:V\ra M$ is lifted to 
$j^*:V[G]\ra M[G^*]$ in $V[G\ast H],$ where $G^*=G\ast H \ast H_o.$
The embedding $j^*$ is obviously again given by an $(\k,\l)$-extender.

To construct a $j^*Q_\k$-generic/$M[G^*]$ filter $H^*\in V[G\ast H]$
consider an enumeration 
$\langle f_\a; \a<\k^+ \rangle$ of all functions from
$\k$ into $[Q_\k]^\k,$ each with cofinally many repetitions.
We need $[H\res\a]^j \in H^*$ for all $\a<\k^+,$ hence
construct a descending sequence 
$\langle p_\a; \a<\k^+ \rangle \subseteq j^*Q_\k$
so that
\vskip5pt

(i) $p_\a \in j^*(Q_\k \res \a),$

(ii) $p_\a \leq [H\res\a]^j,$

(iii) $p_\a$ meets any maximal antichain $\subseteq j^*(Q_\k\res\a)$
of the form $(j^*f_\a)(a)$ for an $a\in [\l]^{<\omega} .$
\vskip5pt

\noindent Since any maximal antichain in $j^*Q_\k$
is actually an antichain in $j^*(Q_\k\res\a)$ for some $\a<\k^+,$
the sequence generates a desired $H^*\in V[G\ast H]$
$j^*Q_\k$-generic over $M[G^*].$
Therefore $j^*$ is lifted to $j^{**}:V[G\ast H]\ra M[G^*\ast H^*].$
Obviously $V[G\ast H] \cap\;^\k M[G^*\ast H^*]\subseteq M[G^*\ast H^*]$
as $P_{\k+1}$ is $\k^+$-c.c. To prove that $j^{**}$ is $\gamma$-strong
it is enough to show that 
$\Cal P^{V[G\ast H]}_\gamma (\k^+) \subseteq M[G^*\ast H^*].$
For each $\delta<\gamma$ fix a bijection
$\pi_\delta:\Cal P_\delta(\k^+)\times P_{\k+1} \ra \Cal P_\delta(\k^+)$
that is in $M$ ($\Cal P_0(\k^+)=\k^+,$ 
$\Cal P_{\delta+1}=\Cal P(\Cal P_\delta)$).
We actually need $\langle \pi_\delta;\delta <\gamma \rangle \in M .$
Then for each element $x$ of $\Cal P^{V[G\ast H]}_\gamma (\k^+)$
use $\pi_\delta$'s to find a code in $\Cal P_\gamma(\k^+)\subseteq M$
for its $P_{\k+1}$-name $\dot x.$
Consequently the name $\dot x$ itself can be decoded in $M,$
and so $x=i_{G\ast H} (\dot x)$ is in $M[G\ast H]\subseteq M[G^*\ast H^*].$

We say that $\k$ is strong if it is $\gamma$-strong for every $\gamma .$
As in the case of a supercompact cardinal we can assume without loss 
of generality that there is no inaccessible cardinal above $\k,$
and then use the same argument to find $P_{\k+1}$ that works
for class many $\gamma$'s preserving the strongness of $\k.$
That concludes our proof of the main theorem.

\enddocument
\end